\documentclass[a4paper,12pt]{article}
\usepackage{amsmath, amssymb,amsfonts}
\title{Galois door Hurwitz uitleg in dagboek 23 }
\author{}
\date{}
\usepackage{amsmath,amssymb}
\usepackage{graphicx}
\usepackage{titlesec}
\titleformat{\section}
{\normalfont\normalsize\bfseries}
{\thesection}{1em}{}

\author{
	Math Dicker, Hoensbroek, the Netherlands 
	}

\date{\today}

\title{Adolf Hurwitz and the Fundamental Theorem of Galois Theory: The Königsberg Lectures of 1890--1891}

\begin{document}
		
		\maketitle

		\begin{abstract}
			In the winter semester of 1890--1891 Adolf Hurwitz delivered a lecture course
			at the Albertina University in Königsberg entitled
			\emph{Theorie der algebraischen Gleichungen}.
			These lectures contain a particularly clear presentation of the ideas of
			Évariste Galois and, in particular, a proof of the fundamental theorem of
			Galois theory formulated in the language of substitutions.
			The present paper analyzes Hurwitz’s treatment of this result on the basis of his lecture notes preserved in the ETH Library in Zurich (Hs~582:66), together with material from his \emph{Mathematisches Tagebuch~23} (Hs~582:23).
			After placing the Königsberg lectures in their historical context,
			we give an overview of their mathematical content and reconstruct in detail
			Hurwitz’s argument leading to the fundamental theorem.
		\end{abstract}

	\section{Historical Context and Sources}
	 
	It seems reasonable to assume that someone who teaches Galois theory has studied Évariste Galois’ \emph{Mémoire sur les conditions de résolubilité des équations par radicaux}, or has at least attempted to do so. Anyone presenting a theory generally benefits from being aware of its original formulation and argumentation. It is well known that the \emph{Mémoire} is difficult to access; in what follows I will show how Adolf Hurwitz (1859–1919) can help to make it more accessible. \
	
	Recently, a series of lecture notes by the student Emil Leutenegger (1894–1978) became available. Between 1915 and 1917 he attended lectures in Zürich given by Adolf Hurwitz. These are extensive notes from the courses \emph{Algebraischen Gleichungen}, \emph{Zahlentheorie}, and \emph{Funktionentheorie}. Substantial material is known for the courses \emph{Zahlentheorie} and \emph{Funktionentheorie}. Numerous notes by Hurwitz himself exist, and both courses have been published in book form: the first by Nikolaos Kritikos, the second by Richard Courant. \
	
	In this article, we focus  on the course \emph{Algebraischen Gleichungen}. Thanks to the efforts of Mrs. Evelyn Boesch (ETH Library, University Archives of ETH Zürich), the surviving manuscripts of Hurwitz for this course have been digitized and made accessible as PDFs. Among this material is a lecture course from the winter semester 1890--1891 in Königsberg \cite{Hurwitz1}.	 \
	
	The course \emph{Theorie der algebraischen Gleichungen} \cite{Hurwitz1} in Königsberg contained the following chapters: \emph{The algebraic solution of equations of the second, third and fourth degree; the theory of symmetric functions; basic concepts of substitution theory; the insolvability of equations of the fifth degree; Galois theory; application of Galois theory to the algebraic solution of equations.} 	\

	In the lectures in Zürich the section on Galois theory is absent, except for the course of 1909      \cite[Chap.~9--16]{Hurwitz2}. Around that time a note also appears in his \emph{Mathematisches Tagebuch No.23} \cite[p.~153]{HurwitzDiary}, in which he sets out the main theorem of Galois theory; a marginal note referring to this is also recorded in the manuscript of the Königsberg lectures \cite[p.~37]{Hurwitz1}. It appears that in 1909 he consulted the Königsberg lecture notes. We focus on that section of his \emph{Mathematisches Tagebuch No.23}. A photograph of this section is included.\
	
	The lecture courses in Königsberg (1890–1891) and Zürich (1909) occupy a special place in the teaching of Galois theory, particularly because they were delivered by Adolf Hurwitz. He still treats the material entirely within the classical nineteenth-century tradition of substitutions. In this way he connects directly to the original \emph{Mémoire} of Évariste Galois, whose central ideas he carefully follows and explicates.\
	
	In doing so he assumes that Galois—without explicitly formulating it— already possessed a notion of what was later called a \emph{Rationalitätsbereich} and by Dedekind a \emph{field} (\emph{Körper}). He consistently refers his students to the contemporary literature, including in this case a German translation of Galois’ \emph{Mémoire} and the book by Eugen Netto, \emph{Substitutionentheorie und ihre Anwendungen auf Algebra} \cite{netto}.\
		
	The notes of this course are preserved in the ETH Library Zürich (Hs 582:66). For the work of Galois, Hurwitz refers to \emph{Abhandlungen über die algebraische Auflösung der Gleichungen} by N.H.~Abel and E.~Galois, which were republished in German translation in 1889 by H.~Maser \cite{maser}. 
		
	In the present article we work out the proof of the main theorem of Galois theory as Hurwitz recorded it in his diary. Although the proof is rather subtle, it is at the same time of great beauty. We prove the theorem in a concrete situation so that the complexity is somewhat reduced, while the essence of the proof is nevertheless preserved.

		\section{Reconstruction of Hurwitz's Proof}

		\begin{figure}[htbp]
			\centering
			\includegraphics[width=1.25\textwidth]{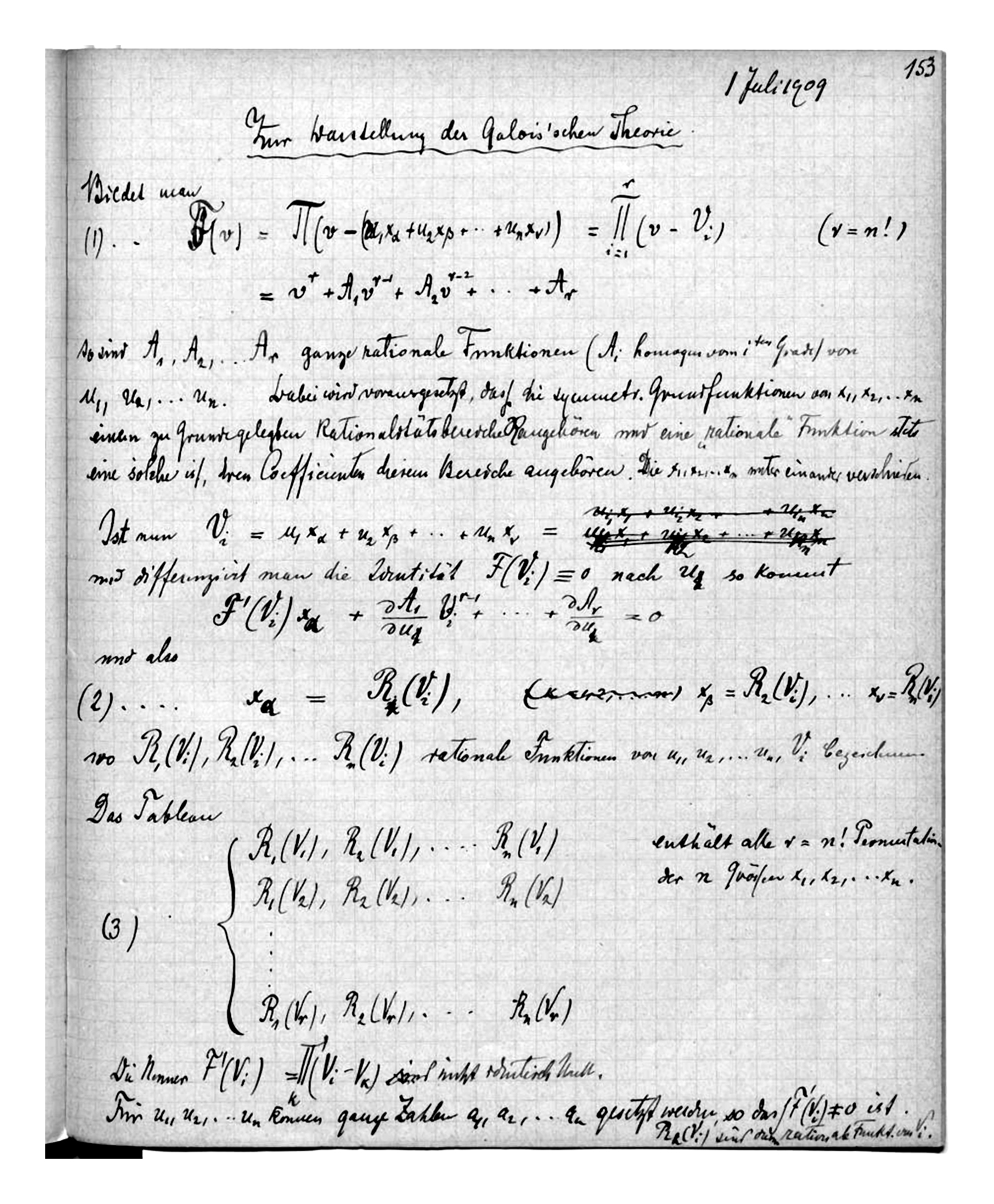}
			\caption{Hurwitz, Mathematisches Tagebuch 23, p.~153; A precise elaboration of a passage from the Königsberg lecture related to Proposition I of Galois.}
			\label{fig:hurwitz1}
		\end{figure}
		
		\begin{figure}[p]
			\centering
			\includegraphics[width=1.25\textwidth]{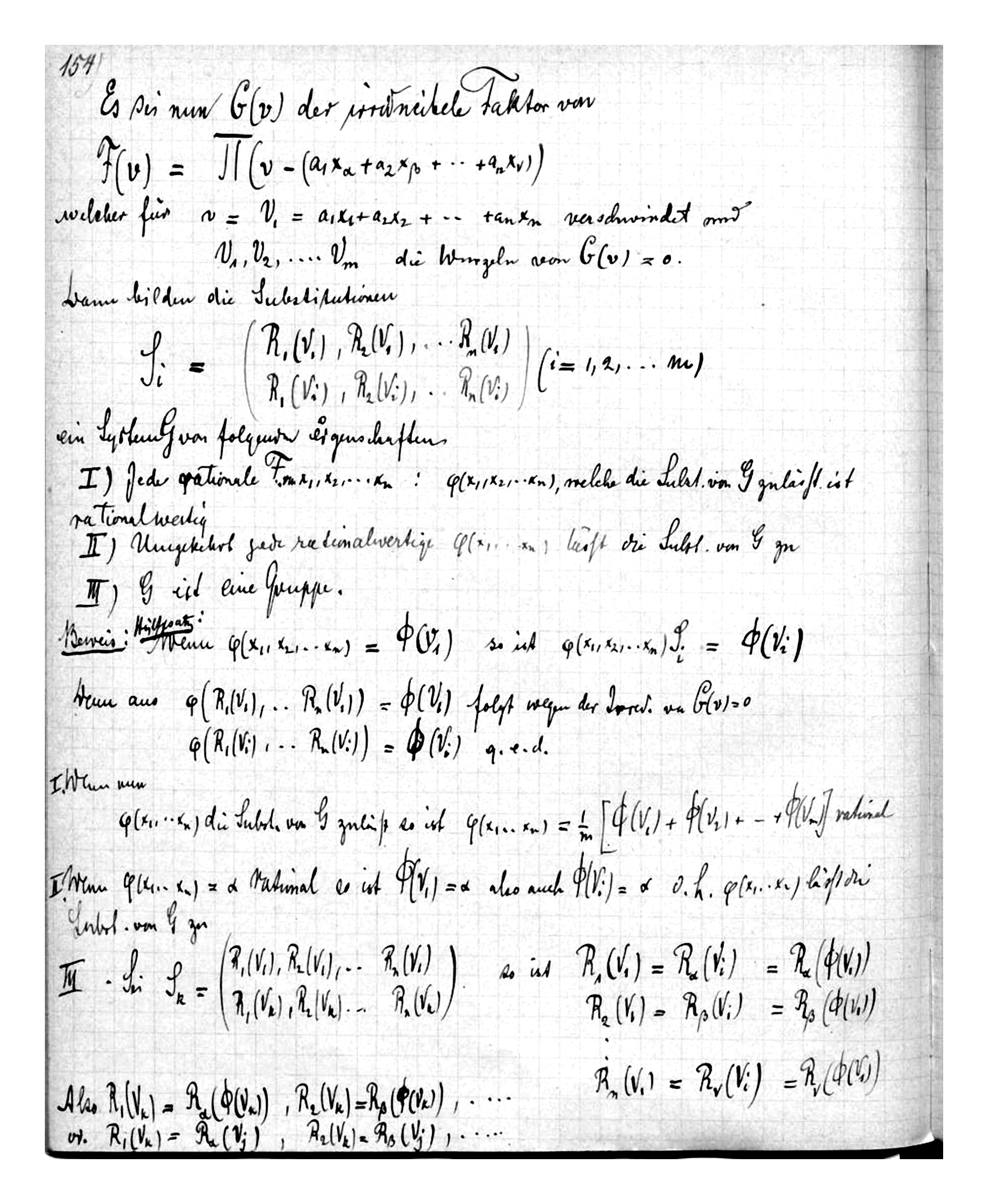}
			
			\vspace{1em}
			
			\includegraphics[width=1.25\textwidth]{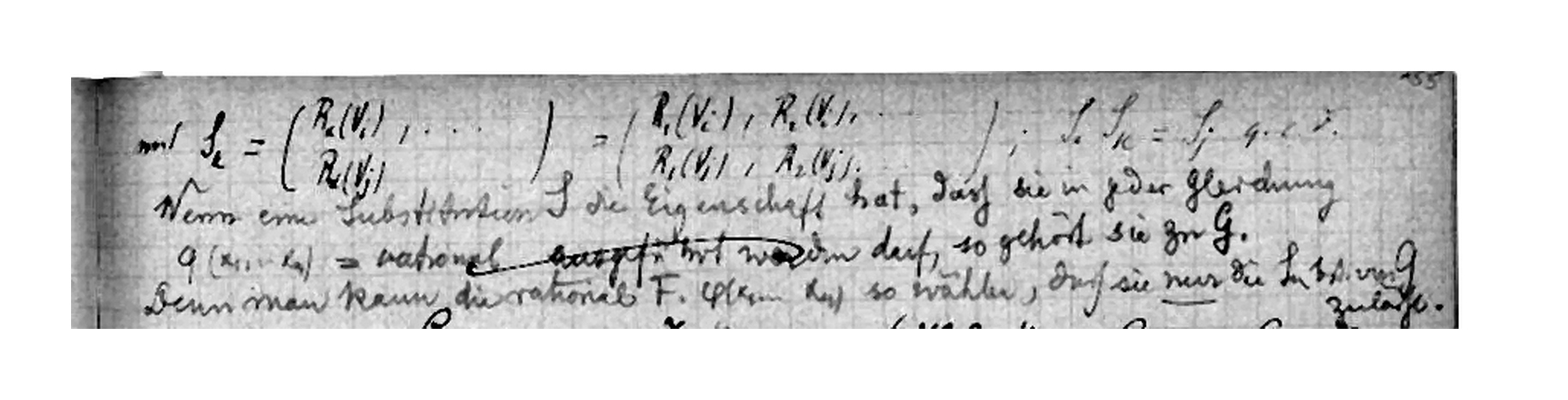}
			
			\caption{Hurwitz, \emph{Mathematisches Tagebuch} 23, pp.~154–155.}
			\label{fig:hurwitz23}
		\end{figure}
		
		The passage from Hurwitz’s diary dated 1 July 1909, in which he proves the main theorem, is reproduced below. First we sketch the context in which the discussion is placed; for convenience we take the base field to be $\mathbb{Q}$.  
		Let $f(x)\in\mathbb{Q}[x]$ be a polynomial with four distinct roots and let  
		$L=\mathbb{Q}(x_1,x_2,x_3,x_4)$ be the splitting field of $f$.  
		
		According to Lemma~II in the \emph{Mémoire} of Galois there exist natural numbers  
		$n_1,n_2,n_3,n_4$ such that the 24 values		
		\(
		n_1x_{s(1)}+n_2x_{s(2)}+n_3x_{s(3)}+n_4x_{s(4)}
		\)		are pairwise distinct, where $s$ permutes the numbers $1,2,3,4$, starting from
				\(
		V_1=n_1x_1+n_2x_2+n_3x_3+n_4x_4 .
		\)
		All these values belong to $L$. We denote them by
		\[
		V_1,V_2,V_3,\dots,V_{24}.
		\]
		
		Given a value $V_i$, one may associate with it a substitution from $S_4$.  
		For example,
		\(
		n_1x_3+n_2x_1+n_3x_4+n_4x_2
		\)
				corresponds to the substitution
		
		\[
		S=
		\begin{pmatrix}
			x_1 & x_2 & x_3 & x_4 \\
			x_3 & x_1 & x_4 & x_2
		\end{pmatrix}.
		\]
		
		We now assume a field $\mathfrak{R}$ contained in $L$.  
		The fact that $f\in\mathfrak{R}[x]$ will be important later in connection with the application of Lemma~I of Galois.
		
		Our aim is the following.  
		The elements
		\(
		V_1,V_2,V_3,\dots,V_{24}
		\)
		are numbers in $L$.  
		We wish to determine functions		
		\(
		R_1,R_2,R_3,R_4\in\mathbb{Q}(x)
		\)
		which yield the values $x_1,x_2,x_3,x_4$ respectively, such that for every $i=1,\dots,24$
		
		\[
		V_i=n_1R_1(V_i)+n_2R_2(V_i)+n_3R_3(V_i)+n_4R_4(V_i).
		\]
		
		The purpose is therefore to reconstruct from the numerical value $V_i$ the original linear combination
		\(
		n_1R_1(V_i)+n_2R_2(V_i)+n_3R_3(V_i)+n_4R_4(V_i)
		\)
				and thereby also the corresponding substitution $S_i$
		\[
		S_i=
		\begin{pmatrix}
			R_1(V_1), & R_2(V_1), & R_3(V_1), & R_4(V_1)\\
			R_1(V_i), & R_2(V_i), & R_3(V_i), & R_4(V_i)
		\end{pmatrix}
		\qquad (i=1,2,\dots,24).
		\]
		
		This line of reasoning appears in Lemma~II, Lemma~III and Proposition~I of Galois’ \emph{Mémoire}.  
		Below follows an adapted transcription of the relevant part of Hurwitz’s notebook.
		
		\newpage
		
		\begin{flushright}
			1~July~1909 \hfill p.~153
		\end{flushright}
		
		\begin{center}
			\textbf{On the treatment of Galois theory}
		\end{center}
		
		Define
		
		\[
		F(v)
		=\prod\bigl(v-(u_1x_{\alpha}+u_2x_{\beta}+u_3x_{\gamma}+u_4x_{\delta})\bigr)
		=\prod_{i=1}^{24}(v-V_i),
		\qquad (r=24)
		\]
		
		\[
		= v^{r}+A_1v^{r-1}+A_2v^{r-2}+\dots +A_r .
		\]
		
		The coefficients \(A_1,\dots,A_r\) are rational functions in the variables
		\(u_1,\dots,u_4\); moreover \(A_i\) is homogeneous of degree \(i\).
		
		We assume that the elementary symmetric functions of
		\(x_1,\dots,x_4\) belong to the field \(\mathfrak{R}\).  
		By a rational function we always mean a function whose coefficients lie in \(\mathfrak{R}\).  
		The elements \(x_1,\dots,x_4\) are mutually distinct.
		
		Let
		
		\[
		V_i=u_1x_\alpha+u_2x_\beta+u_3x_\gamma+u_4x_\delta .
		\]
		
		Differentiating the identity \(F(V_i)\equiv0\) with respect to \(u_1\) yields
		
		\[
		F'(V_i)\cdot x_\alpha
		+\frac{\partial A_1}{\partial u_1}V_i^{r-1}
		+\dots
		+\frac{\partial A_r}{\partial u_1}
		=0 .
		\]
		
		From this one obtains analogously
		
		\[
		x_\alpha=R_1(V_i),\qquad
		x_\beta=R_2(V_i),\dots,x_\delta=R_4(V_i).
		\]
		
		The scheme
		
		\[
		\begin{cases}
			R_1(V_1),\dots,R_4(V_1)\\
			R_1(V_2),\dots,R_4(V_2)\\
			\vdots\\
			R_1(V_r),\dots,R_4(V_r)
		\end{cases}
		\]
		
		contains all \(r=24\) permutations of the four elements
		\(x_1,\dots,x_4\).
		
		The denominators
		
		\[
		F'(V_i)=\prod_{k\neq i}(V_i-V_k)
		\]
		
		are not identically zero.
		
		One can choose integers \(a_1,\dots,a_n\) such that
		\(F'(V_i)\neq0\).  
		In that case the \(R_k(V_i)\) are rational functions of \(V_i\).
		
		\begin{flushright}
		p.~154
		\end{flushright}
		
		Let \(G(v)\) be the irreducible factor of
		
		\[
		F(v)
		=\prod\bigl(v-(a_1x_\alpha+a_2x_\beta+a_3x_\gamma+a_4x_\delta)\bigr),
		\]
		
		which vanishes for
		
		\[
		v=V_1=a_1x_1+a_2x_2+a_3x_3+a_4x_4,
		\]
		
		and whose zeros (after possible renumbering) are
		
		\[
		V_1,\dots,V_m .
		\]
		
		Then the \(m\) substitutions
		
		\[
		S_i=
		\begin{pmatrix}
			R_1(V_1) & \dots & R_4(V_1)\\
			R_1(V_i) & \dots & R_4(V_i)
		\end{pmatrix}
		\qquad (i=1,\dots,m)
		\]
		
		define a system \(G\) with the following properties:
		
		\begin{description}
			
			\item[I)]  
			Every rational function \(F\) of \(x_1,x_2,x_3,x_4\) whose value
			\(\varphi(x_1,x_2,x_3,x_4)\)
			admits the substitutions of \(G\) has a value belonging to \(\mathfrak{R}\).
			
			\item[II)]  
			Every rational function \(F\) of \(x_1,x_2,x_3,x_4\) whose value
			\(\varphi(x_1,x_2,x_3,x_4)\)
			belongs to \(\mathfrak{R}\) admits the substitutions of \(G\).
			
			\item[III)]  
			\(G\) is a group.
			
			\item[IV)]  
			A substitution \(S\) belongs to \(G\) if and only if every rational function
			\(\varphi(x_1,x_2,x_3,x_4)\) belonging to \(\mathfrak{R}\)
			admits the substitution \(S\).
			
		\end{description}
		
		Hurwitz proves these statements using an auxiliary lemma and shows that
		the correspondence between rationality domains and substitution groups
		yields precisely the fundamental theorem of Galois theory.\\

		Proof. 
		
			\begin{description}
		\item[\textbf{Auxiliary Lemma.}] 
		Assume that
		\(
		\varphi(x_1,x_2,\dots,x_n)=\Phi(V_1)	\footnote{
			Since \(x_i=R_i(V_1)\) and the functions \(R_i\) belong to the rationality domain \(\mathbb{Q}(x)\), we obtain
			\(
			\mathbb{Q}(x_1,x_2,\ldots,x_n)=\mathbb{Q}(V_1).
			\)
			Hence \(\Phi\) may be chosen as a polynomial in \(\mathbb{Q}[x]\). 
		}
		\)
		then	\(
		S_i(\varphi(x_1,x_2,\dots,x_n))=\Phi(V_i).
		\)
			Indeed,		\(
		\varphi\bigl(R_1(V_1),\dots,R_n(V_1)\bigr)=\Phi(V_1)
		\)		implies, by the irreducibility of \(G(v)\), that		\(
		\varphi\bigl(R_1(V_i),\dots,R_n(V_i)\bigr)=\Phi(V_i).
		\)

		\item[\textbf{I)}]
		
		Suppose that the function
		\(
		\varphi(x_1,\dots,x_n)
		\)
		is invariant under the substitutions of \(G\). Then
				\(
		\varphi(x_1,\dots,x_n)
		=
		\frac{1}{m}\bigl[\Phi(V_1)+\Phi(V_2)+\dots+\Phi(V_m)\bigr]
		\)
		
		is rational\footnote{
			Since \(\Phi \in \mathfrak{R}[x]\) and the Newton identities for
			\(V_1,\dots,V_m\) can be expressed in terms of the coefficients of
			\(G(v)\), the expression above lies in \(\mathfrak{R}\).
		}.

		\item[\textbf{II)}]
		
		If		
		\(
		\varphi(x_1,\dots,x_n)=\alpha
		\)
		is rational, then \(\Phi(V_1)=\alpha\) and therefore also \(\Phi(V_i)=\alpha\). 
		Consequently \(\varphi(x_1,\dots,x_n)\) is invariant under the substitutions of \(G\).
		
		\item[\textbf{III)}]
		
		We show that
				\(
		S_k\circ S_i = S_j ,
		\)
				where \(1\le i,k,j\le m\) and \(m\) denotes the degree of \(G(v)\).
				Consider
				\[
		S_k=
		\begin{pmatrix}
			R_1(V_1),R_2(V_1),R_3(V_1),R_4(V_1)\\
			R_1(V_k),R_2(V_k),R_3(V_k),R_4(V_k)
		\end{pmatrix}.
		\]
		
		Since \(V_i=\Phi(V_1)\) with \(\Phi(x)\in\mathbb{Q}[x]\), we obtain
			\[
			R_1(V_1)=R_\alpha(V_i))=R_\alpha(\Phi(V_1))\quad
			R_2(V_1)=R_\beta(V_i))=R_\beta(\Phi(V_1))\quad\dots
			\]
			
		and subsequently	
		\[
		R_1(V_k)=R_\alpha(\Phi(V_k)),\quad
		R_2(V_k)=R_\beta(\Phi(V_k)),\dots
		\]
		
		for suitable indices \(\alpha,\beta,\gamma,\delta\).
			Applying the auxiliary lemma to \(V_i=\Phi(V_1)\), we obtain
			\(S_k(V_i)=\Phi(V_k)\).	
		Moreover \(\Phi(V_k)\) is a root of \(G(v)\), since
				\(
		G(\Phi(V_1))=0
		\)
				implies
				\(
		G(\Phi(V_k))=0.
		\)
			Hence \(S_k(V_i)\) is one of the roots
				\(
		V_1,V_2,\dots,V_m.
		\)	Thus there exists \(j\) such that
				\(
		S_k(V_i)=V_j .
		\)
				Consequently
				\(
		R_1(V_k)=R_\alpha(V_j),\quad
		R_2(V_k)=R_\beta(V_j),\dots
		\)
			
		\[
	S_k=
	\begin{pmatrix}
		R_\alpha(V_i),R_\beta(V_i),R_\gamma(V_i),R_\delta(V_i)\\
		R_\alpha(V_j),R_\beta(V_j),R_\gamma(V_j),R_\delta(V_j)
	\end{pmatrix}
	\]
	
		\[
	S_k=
	\begin{pmatrix}
		R_1(V_i),R_2(V_i),R_3(V_i),R_4(V_i)\\
		R_1(V_j),R_2(V_j),R_3(V_j),R_4(V_j)
	\end{pmatrix}
	\]

				and therefore
			\(
		S_k \circ S_i = S_j .
		\)
		
		\item[\textbf{IV)}]
		
		Let \(S\) be a substitution that can be applied to every equation of the form
				\(
		\varphi(x_1,\ldots,x_n)=\text{rational},
		\)
			then \(S\) belongs to \(G\). Indeed, one can choose a rational function
		\(\varphi(x_1,\ldots,x_n)\)		\footnote{ \cite[Chap.~10, \S 4]{Hurwitz2} Let $T$ be a group of substitutions in $S_4$ and let
			\(
			V_1 = n_1x_1 + n_2x_2 + n_3x_3 + n_4x_4
			\)
			be chosen such that the $24$ values obtained from $V_1$ by applying all substitutions of $S_4$ are pairwise distinct.		Consider the polynomial
			\(
			\psi(t)=\prod_{\tau\in T}(t-\tau(V_1))
			\)
			and for $x\in S_4$ define
			\(
			\psi_x(t)=\prod_{\tau\in T}(t-(x\circ\tau)(V_1)).
			\)
			Then
			\(
			\psi_x(t)=\psi(t) \iff x\in T .
			\)
						Let
			\(
			\theta_x(t)=\psi_x(t)-\psi(t).
			\)
			For every $x\in S_4$ with $x\notin T$, the polynomial $\theta_x(t)$ is nonzero and therefore has only finitely many zeros. 
			Hence there exists a natural number $n$ such that
			\(
			\theta_x(n)\neq 0 \qquad \text{for all } x\notin T .
			\)
						Put
			\(
			\alpha=\psi(n),
			\)
		then
			\(
			x(\alpha)=\alpha \iff x\in T .
			\)}
				 that admits only the substitutions belonging to \(G\)
				\footnote{
			This corresponds to the classical correspondence between subgroups and rationality domains. 
			Hurwitz associates with each rationality domain \(\mathfrak{R}\subset K(x_1,\dots,x_n)\) a uniquely determined group \(G\). 
			This group satisfies \(L^G=\mathfrak{R}\). In particular, \(G\) is the Galois group of \(K(x_1,\dots,x_n)\) over \(\mathfrak{R}\).
		}

			\section{Concluding Remarks}

	During the lecture in 1909 Hurwitz addressed his students as follows:
	
	\textit{\textbf{Galois ist einer der größten mathematischen Genie's aller Zeiten.}
		Geboren 1811 gestorben in Folge eines Duells 1832 ist er nur
		20$\tfrac{1}{2}$ Jahr alt geworden. Schon mit 19 Jahren war er in Besitz
		der nach ihm benannten Theorie, welche der Lehre von den algebraischen
		Gleichungen ein ganz neues Gepräge verlieh.}\\
		
	Two lecture courses by Hurwitz devoted to Galois theory are known.
	The first was given in Königsberg during the winter semester of 1890--1891, when Hurwitz was about thirty-one years old; the second was delivered in Zürich in the summer semester of 1909, when he was about fifty.
	The Königsberg course is the more extensive one and follows quite closely Galois’ difficult Mémoire.
	
	In 1909 Hurwitz returned to the Königsberg material. He recorded this proof separately in his diary from that period and added a corresponding marginal note in the manuscript of the Königsberg lectures.
	
	The proof is subtle, concise, and elegant; it is reproduced in the present article with only minor adaptations.
	The Königsberg lectures are of particular interest, not only because they offer insight into Hurwitz as a teacher—indeed an excellent one—but also because they provide a valuable point of access to \emph{Galois’ Mémoire}. In this way the Königsberg lectures form an illuminating bridge between Galois’ original ideas and their later exposition in Hurwitz’s teaching.
			
		\end{description}

\end{document}